# REFLEXIVE SUBGROUPS OF THE BAER-SPECKER GROUP AND MARTIN'S AXIOM

RÜDIGER GÖBEL AND SAHARON SHELAH

ABSTRACT. In two recent papers [9, 10] we answered a question raised in the book by Eklof and Mekler [7, p. 455, Problem 12] under the set theoretical hypothesis of $\diamondsuit_{\aleph_1}$ which holds in many models of set theory, respectively of the special continuum hypothesis (CH). The objects are reflexive modules over countable principal ideal domains $R$, which are not fields. Following H. Bass [1] an $R$-module $G$ is reflexive if the evaluation map $\sigma : G \longrightarrow G^{**}$ is an isomorphism. Here $G^* = \text{Hom}(G, R)$ denotes the dual module of $G$. We proved the existence of reflexive $R$-modules $G$ of infinite rank with $G \not\cong G \oplus R$, which provide (even essentially indecomposable) counter examples to the question [7, p. 455]. Is CH a necessary condition to find 'nasty' reflexive modules? In the last part of this paper we will show (assuming the existence of supercompact cardinals) that large reflexive modules always have large summands. So at least being essentially indecomposable needs an additional set theoretic assumption. However the assumption need not be CH as shown in the first part of this paper. We will use Martin's axiom to find reflexive modules with the above decomposition which are submodules of the Baer-Specker module $R^\omega$.

## 1. INTRODUCTION

We will derive our results for abelian groups, but it is an easy exercise to replace the ground ring $\mathbb{Z}$ by any countable principal ideal domain which is not a field. Just notice that we could work with one prime only! For supercompact cardinals we refer either to Jech [13] or to Kanamori [14]. If $G$ is any abelian group then $G^* = \text{Hom}(G, \mathbb{Z})$ denotes its dual group, and $G$ is a dual if $G \cong D^*$ for some abelian group $D$.

Particular dual groups are the reflexive groups $D$, see Bass [1, p. 476]. Recall that
$$\sigma = \sigma_D : D \longrightarrow D^{**} \ \ (d \longrightarrow \sigma(d))$$
with $\sigma(d) \in D^{**}$ and
$$\sigma(d) : D^* \longrightarrow \mathbb{Z} \ \ (\varphi \longrightarrow \varphi(d))$$
is the evaluation map and $D$ is *reflexive* if the evaluation map $\sigma_D$ is an isomorphism. Recent results about reflexive and dual abelian groups are discussed in [7, 9, 10].

1991 *Mathematics Subject Classification.* primary: 13C05, 13C10, 13C13, 20K15, 20K25, 20K30; secondary: 03E05, 03E35.

This work is supported by the project No. G-545-173.06/97 of the German-Israeli Foundation for Scientific Research & Development
GbSh 727 in Shelah's list of publications.

In the third section we will show that dual groups, in particular reflexive groups may have large summands, hence can't be essentially indecomposable without any set-theoretic restrictions.

**Theorem 1.1.** *If $\kappa$ is a supercompact cardinal and $H$ is a dual group of cardinality $\geq \kappa$, then there is a direct summand $H'$ of $H$ with $\chi \leq |H'| < \kappa$ for any cardinal $\chi < \kappa$.*

This theorem explains that we had to use CH in [9] and alternatively below we will use Martin's axiom (and possibly negation of CH).

In order to prove a result in contrast of Theorem 1.1 we use scalar products on the Baer-Specker group $P$. Recall that

$$P = \mathbb{Z}^\omega$$

is the set of all elements

$$\mathbf{x} = \sum_{i \in \omega} x_i \mathbf{e}_i \text{ with } x_i \in \mathbb{Z}$$

and $\mathbf{e}_i \in P$ defined by the Kronecker symbol, addition is defined component-wise. Throughout this paper we will adopt the convention in writing elements of $P$ as displayed in the last formula. The Baer-Specker group $P$ has the subgroup $S$ of all elements $\mathbf{x}$ of finite support, that is $x_i = 0$ for almost all $i \in \omega$. The crucial subgroup for constructing reflexive groups is the $\mathbb{Z}$-adic closure $\mathbb{D}$ of $S$ in $P$. This will be our target in Section 3. We will also show that the the endomorphism ring of such a reflexive abelian group can be $\mathbb{Z}$ *modulo* the ideal of all endomorphisms of finite rank. We have the the following

**Theorem 1.2.** *(ZFC + MA) There are two subgroups $H_i$ ($i = 1, 2$) of the Baer-Specker group $P$ with the following properties.*
  (i) *$S \subseteq H_i \subseteq_* \mathbb{D}$ are pure.*
 (ii) *$H_i$ is $\aleph_1$-free and slender.*
(iii) *There is a natural bilinear form $\Phi : H_1 \times H_2 \longrightarrow \mathbb{Z}$ arriving from scalar product on $P$ which induces $H_1^* \cong H_2$ and $H_2^* \cong H_1$ such that $H_1, H_2$ are reflexive.*
 (iv) *$H_i \oplus \mathbb{Z} \not\cong H_i$ for $i = 1, 2$.*
  (v) *$\text{End } H_i = \mathbb{Z} \oplus \text{Fin } H_i$.*

Note that $\varphi \in H_1^*$ is induced by $\Phi$ if there is $h \in H_2$ such that $\varphi = \Phi(\ , h)$. The set $\text{Fin } H_i$ of all endomorphisms of $H_i$ with finite rank image is an ideal of the endomorphism ring $\text{End } H_i$ and the last statement of the theorem means that this ideal is a split extension in $\text{End } H_i$.

Hence each $H_i$ is separable and essentially indecomposable, this means any decomposition $H_i = C \oplus E$ must have a summand $E$ or $C$ of finite rank. The key for proving Theorem 1.2 are new algebraic and combinatorial methods and some old techniques from earlier papers like [11] or [4].

## 2. Reflexive groups of cardinality $\leq 2^{\aleph_0}$ under Martin's axiom

In contrast to the results in Section 3 concerning the existence of arbitrarily large summands of reflexive groups larger then a supercompact cardinality in this section we will now construct essentially indecomposable reflexive groups under Martin's axiom MA. As above let $P = \prod_{n \in \omega} \mathbf{e}_n \mathbb{Z}$ be the Baer-Specker group of all elements

$$P = \{\mathbf{x} = \sum_{i \in \omega} x_i \mathbf{e}_i : \ (x_i \in \mathbb{Z})\}.$$

Here $\mathbf{e}_i$ can be viewed as the element $\mathbf{x}$ with coefficients $x_{ij} = \delta_{ij}$ the Kronecker symbol. Hence

$$S = \langle \mathbf{e}_i : i \in \omega \rangle = \bigoplus_{i \in \omega} \mathbf{e}_i \mathbb{Z}$$

is a subgroup of $P$ of all elements $\mathbf{x}$ of finite support

$$[\mathbf{x}] = \{i \in \omega : x_i \neq 0\}$$

and $P/S$ is algebraically compact by an old result of Balcerzyk, see Fuchs [8]. Obviously $P/S$ is torsion-free or equivalently $S$ is pure in $P$. Pure subgroups $X \subseteq P$ are denoted by $X \subseteq_* P$. Moreover, let $\mathbb{D}$ be the $\mathbb{Z}$-adic closure of $S$ in $P$, hence $\mathbb{D}/S$ is the maximal divisible (torsion-free) subgroup of $P/S$ which has size $2^{\aleph_0}$. If $H$ is an abelian group, then $\mathrm{Fin}\, H$ denotes the ideal of all endomorphisms $\sigma \in \mathrm{End}\, H$ with $\mathrm{Im}\, \sigma$ of finite rank. The groups we want to construct will be sandwiched between $S$ and $\mathbb{D}$.

We will use Martin's axiom for $\sigma$-centered sets, which is a (proper) consequence of the well-known Martin's axiom and equivalent to the combinatorial principle $P(2^{\aleph_0})$ (below) as shown by Bell [2]. Recall that $D \subseteq \mathfrak{P}$ is dense in the poset $\mathfrak{P}$ if for any $p \in \mathfrak{P}$ there exists $d \in D$ such that $p \leq d$. Martin's axiom is based on posets $\mathfrak{P}$ with c.c.c. using that $p, q \in \mathfrak{P}$ are compatible if there is $r \in \mathfrak{P}$ with $\{p, q\} \leq r$. Recall that $F \subseteq \mathfrak{P}$ is bounded by $r$, say $F \leq r$ if $f \leq r$ for all $r \in F$. A set $X \subseteq \mathfrak{P}$ is directed if all finite subsets of $X$ are bounded in $X$ and $X$ is called $\sigma$-centered ($\sigma$-directed) if it is the countable union of directed subsets. Replacing c.c.c. by '$\sigma$-centered' MA turns into *Martin's axiom for $\sigma$-centered sets:*

Let $\mathfrak{D}$ be a collection of dense subsets $D$ of the poset $\mathfrak{P}$. If $|\mathfrak{D}| < 2^{\aleph_0}$ and $(\mathfrak{P}, \leq)$ is a $\sigma$-centered poset then there is a $\mathfrak{D}$-generic subset $G \subset \mathfrak{P}$. Hence $G$ is directed and meets every $D \in \mathfrak{D}: \ G \cap D \neq \emptyset$.

See [7, p. 164] for MA with c.c.c. Note that the main result in Bell [2] is that Martin's axiom for $\sigma$-centered sets is equivalent to

The combinatorial principle $P(2^{\aleph_0})$: If $\mathfrak{D}$ is a collection of subsets of $\omega$ such that $|\mathfrak{D}| < 2^{\aleph_0}$ and $\bigcap F$ is infinite for every finite $F \subseteq \mathfrak{D}$, then there is an infinite $B \subseteq \omega$ such that $B \setminus D$ is infinite for all $D \in \mathfrak{D}$.

Martin's axiom will help us to define a scalar product or bilinear form $\Phi$ on suitable pairs $\mathbb{H} = (H_1, H_2)$ of pure subgroups $H_j$ of $\mathbb{D}$. We begin with

$$\Phi : S \times S \longrightarrow \mathbb{Z} \text{ with } \Phi(\mathbf{e}_i, \mathbf{e}_j) = \delta_{ij}.$$

Hence $\Phi$ is the unique integer valued, bilinear form on $S \times S$. It extends by continuity uniquely to the non-degenerated, symmetric bilinear form

$$\Phi : \mathbb{D} \times \mathbb{D} \longrightarrow \widehat{\mathbb{Z}} \text{ with } \widehat{\mathbb{Z}} \text{ the } \mathbb{Z}\text{-adic completion of } \mathbb{Z}.$$

We keep this map fixed though out this section and also denote restrictions to pairs of subgroups by $\Phi$. Note that $\widehat{\mathbb{Z}}$ is the cartesian product of the additive groups of $p$-adic integers over all primes $p$ and if

$$\mathbf{a} = \sum_{i \in \omega} a_i \mathbf{e}_i \in \mathbb{D} \text{ and } \mathbf{b} = \sum_{i \in \omega} b_i \mathbf{e}_i \in \mathbb{D}, \text{ then } \Phi(\mathbf{a}, \mathbf{b}) = \sum_{i \in \omega} a_i b_i$$

is well-defined and symmetry $\Phi(\mathbf{a}, \mathbf{b}) = \Phi(\mathbf{b}, \mathbf{a})$ is obvious. Now we consider such pairs $\mathbb{H} = (H_1, H_2)$ such that $\Phi \upharpoonright (H_1, H_2)$ takes only values in $\mathbb{Z}$. More precisely, let $\mathbb{H} \in \mathfrak{P}$ if and only if the followings holds for $j = 1, 2$:

(i) $S \subseteq H_j \subseteq_* \mathbb{D}$
(ii) $|H_j| < 2^{\aleph_0}$
(iii) $\Phi : H_1 \times H_2 \longrightarrow \mathbb{Z}$.

We define a partial order on $\mathfrak{P}$ and say

**Definition 2.1.** *If $\mathbb{H}, \mathbb{H}' \in \mathfrak{P}$ then $\mathbb{H} \subseteq \mathbb{H}'$ if and only if $H_1 \subseteq H_1'$ and $H_2 \subseteq H_2'$.*

The next crucial lemma of this paper will show under MA that $\mathfrak{P}$ is a rich structure.

**Main Lemma 2.2.** *(ZFC + MA) Let $\mathbb{H} = (H_1, H_2) \in \mathfrak{P}, \mathbf{b} \in P \setminus \mathbb{D}$ and $\mathbf{b}^n \in H_1$ for $n \in \omega$. Then there is $\mathbf{a} = \sum_{i \in \omega} a_i \mathbf{e}_i \in \mathbb{D}$ such that for $H_1' = \langle H_1, \mathbf{a} \rangle_* \subseteq \mathbb{D}$ and $\mathbb{H}' = (H_1', H_2)$ the following holds.*

(i) $\mathbb{H} \subseteq \mathbb{H}' \in \mathfrak{P}$ *and* $\Phi(\mathbf{a}, \mathbf{b}) \in \widehat{\mathbb{Z}} \setminus \mathbb{Z}$.
(ii) (a) *Either* $\sum_{i \in \omega} a_i \mathbf{b}^i \notin H_1'$
  (b) *or there is $t \in \mathbb{Z}$ such that $\langle \mathbf{b}^j - t \mathbf{e}_j : j \in \omega \rangle$ is a free direct summand of finite rank.*

**Remark.** By symmetry we obtain a dual result of the Main Lemma 2.2 with $\mathbf{a} \in H_2'$ and $\Phi(\mathbf{b}, \mathbf{a}) \in \widehat{\mathbb{Z}} \setminus \mathbb{Z}$ and $(ii)$ accordingly.

**Proof.** Let $\mathbf{b} = \sum_{i \in \omega} b_i \mathbf{e}_i \in P \setminus \mathbb{D}$ and $\mathbb{H} = (H_1, H_2) \in \mathfrak{P}$ be given by the lemma. Moreover we assume that condition $(ii)(b)$ of the lemma does not hold. This is to say that we must show $(ii)(a)$ of the lemma. This implication will follow at the end of the proof from density of the sets $D^5_{\mathbf{d}tn_0}$ and density will be a consequence of the assumption just made.

We want to approximate $\mathbf{a} \in H_1'$ by a forcing notion $\mathfrak{F}$, a partially ordered set, used for application of MA. The elements $p \in \mathfrak{F}$ are triples

$$(M^p, A^p, n^p) \text{ with } A^p = \langle a_l^p : l < l^p \rangle, \ M^p = \{m_{\mathbf{x}}^p = \sum_{l < l^p} x_l a_l^p : \mathbf{x} = \sum_{i \in \omega} x_i \mathbf{e}_i \in u^p\}$$

subject to the following conditions.

(i) $u^p$ is a finite subset of $H_2$,
(ii) $l^p \in \omega$, $a_l^p, m_{\mathbf{x}}^p \in \mathbb{Z}$, and $n^p \in \mathbb{N}$.

We call $l^p$ the length of the finite sequence of integers $A^p$ and note that $n|m$ means $n$ divides $m$ in $\mathbb{Z}$. In order to turn $\mathfrak{F}$ into a partially ordered set let $p \leq q$ for some $p, q \in \mathfrak{F}$ if the following holds.

$$u^p \subseteq u^q, \ l^p \leq l^q, \ A^p = A^q \upharpoonright l^p,$$

$$n^p | n^q, \text{ and if } l^p \leq l < l^q \text{ then } n^p | a_l^q,$$

if $\mathbf{x} = \sum_{l \in \omega} x_l \mathbf{e}_l \in u^p$ then $m_{\mathbf{x}}^q = m_{\mathbf{x}}^p =: \sum_{l < l^p} x_l a_l^p$ or equivalently $\sum_{l^p \leq l < l^q} x_l a_l^p = 0$.

If $p, q \in \mathfrak{F}$ then let

$$p \sim q \quad \Leftrightarrow \quad (l^p = l^q, A^p = A^q, n^p = n^q)$$

and note that $\sim$ is an equivalence relation on $\mathfrak{F}$. If $p \in \mathfrak{F}$ then let $\mathfrak{F}_p = \{q \in \mathfrak{F} : q \sim p\}$. Surely $\mathfrak{F}$ decomposes into countably many of such uncountable equivalence classes $\mathfrak{F}_p$. We claim that each of them is directed. If $q_1, q_2 \in \mathfrak{F}_p$ then $n^{q_i} = n^p, A^{p_i} = A^p, l^{p_i} = l^p$, hence $q_i = (M^{q_i}, A^p, n^p)$, and if $\mathbf{x} = \sum_{i \in \omega} x_i \mathbf{e}_i \in u^{q_1} \cap u^{q_2}$, then

$$m_\mathbf{x}^{q_1} = \sum_{i < l^{q_1}} x_i a_i^{q_1} = \sum_{i < l^p} x_i a_i^p = m_\mathbf{x}^{q_2}.$$

If we define $q' \in \mathfrak{F}$ by $u^{q'} = u^{q_1} \cup u^{q_2}, A^{q'} = A^p, l^{q'} = l^p, n^{q'} = n^p$, then $M^{q'} = \{m_\mathbf{x}^{q'} = \sum_{i < l^p} x_i a_i^p : \mathbf{x} \in u^{q'}\} = M^{q_1} \cup M^{q_2}$, hence $q' = (M^{q'}, A^p, n^p)$ is a member of $\mathfrak{F}$ and $q_1, q_2 \leq q'$. The claim is shown and by definition

(2.1) $\quad (\mathfrak{F}, \leq)$ is a $\sigma$-centered poset,

as required for applications of MA for $\sigma$-centered sets. In order to apply MA effectively we must define dense subsets of $\mathfrak{F}$ which describe 'local properties' of the desired $\mathbf{a} \in \mathbb{D}$. If $\mathbf{x} = \sum_{i \in \omega} x_i \mathbf{e}_i \in H_2, m \in \mathbb{N}, l_0 \in \omega$ then let

$$D_\mathbf{x}^1 = \{p \in \mathfrak{F} : \mathbf{x} \in u^p\}, \ D_m^2 = \{p \in \mathfrak{F} : m | n^p\}, \ D_{l_0}^3 = \{p \in \mathfrak{F} : l_0 \leq l^p\}$$

,

$$D_m^4 = \{p \in \mathfrak{F} : \sum_{l < l^p} b_l a_l^p \not\equiv m \mod n^p\}$$

and for $\mathbf{d} \in \mathbb{D}, t \in \mathbb{Z}$ and $n_0 \in \mathbb{N}$ let

$$D_{\mathbf{d} t n_0}^5 = \{p \in \mathfrak{F} : \exists m \in \mathbb{N} \ (m | n^p, \ n_0 \sum_{i < l^p} a_i^p \mathbf{b}^i - t \sum_{i < l^p} a_i^p \mathbf{e}_i - \mathbf{d} \not\equiv 0 \mod m\mathbb{D}\}.$$

First note that we defined $< 2^{\aleph_0}$ subsets of $\mathfrak{F}$ as required for MA. Next we want to show that all these sets are dense in $\mathfrak{F}$. The first three cases are easy while the remaining two cases need work. For $D_\mathbf{x}^1$ with $\mathbf{x} = \sum_{i \in \omega} x_i \mathbf{e}_i$ we take any $p \in \mathfrak{F}$ and define $q$ like $p$ just enlarging $u^q = u^p \cup \{\mathbf{x}\}$, let $m_\mathbf{x}^q = \sum_{l < l^p} x_l a_l^p$ and enlarge $M^q = \{m_\mathbf{y}^p : \mathbf{y} = \sum_{i \in \omega} y_i \mathbf{e}_i \in u^p\} \cup \{m_\mathbf{x}^q\}$ as well, hence $q \leq q$ and $D_\mathbf{x}^1$ is dense in $\mathfrak{F}$. Similarly take any $p \leq q \in \mathfrak{F}$ with $m | n^q$, hence $D_m^2$ is dense. For $D_{l_0}^3$ replace any $A^p$ by $A^q = (A^p)^\wedge(0, \ldots, 0)$ with $(0, \ldots, 0)$ a vector of $l_0$ zeros and let $u^q = u^p, l^q = l^p + l_0, n^q = n^p$. In the fourth case we first notice that $\mathbf{b} = \sum_{i \in \omega} b_i \mathbf{e}_i \in P \setminus \mathbb{D}$ by hypothesis, hence there is $s' \in \mathbb{N}$ such that the set

$$W = \{k \in \omega : b_k \in \mathbb{Z} \setminus s'\mathbb{Z}\} \text{ is infinite.}$$

Suppose $p \in \mathfrak{F}$ contradicts the density of $D_m^4$ for some $m \in \mathbb{N}$, hence

(2.2) $\qquad\qquad$ there is no $q \in D_m^4$ with $p \leq q$.

We write

$$u^p = \{\mathbf{a_1}, \ldots \mathbf{a_{k-1}}\} \subseteq H_2 \text{ and let } \mathbf{a_j} = \sum_{i \in \omega} a_{ji} \mathbf{e}_i.$$

Also consider the $k \times \omega$-matrix ($s \in \omega$) (**G**)

$$\begin{pmatrix} a_{11} & a_{12} & \ldots & a_{1s} & \ldots \\ a_{21} & a_{22} & \ldots & a_{2s} & \ldots \\ \vdots & & & & \\ a_{k-1,1} & a_{k-1,2} & \ldots & a_{k-1,s} & \ldots \\ b_1 & b_2 & \ldots & b_s & \ldots \end{pmatrix}$$

as well as the $(k-1) \times \omega$-matrix (**H**) deleting the last row of $b_s$'s of the matrix above.

$$\begin{pmatrix} a_{11} & a_{12} & \ldots & a_{1s} & \ldots \\ a_{21} & a_{22} & \ldots & a_{2s} & \ldots \\ \vdots & & & & \\ a_{k-1,1} & a_{k-1,2} & \ldots & a_{k-1,s} & \ldots \end{pmatrix}$$

We pick finite subsets $w$ of $[l^p, \omega)$ and consider the column vectors $g_l^p$ ($l \in w$) of the first matrix (**G**) and $h_l^p$ ($l \in w$) of the second matrix (**H**) accordingly and claim that for all finite

(2.3) $\quad w \subseteq [l^p, \omega)$ and $d_l \in \mathbb{Q}$ $\quad [\sum_{l \in w} d_l h_l^p = 0 \Leftrightarrow \sum_{l \in w} d_l g_l^p = 0]$.

The proof "$\Leftarrow$" is trivial. For "$\Rightarrow$" suppose for contradiction that

$$\sum_{l \in w} d_l h_l^p = 0 \text{ but } \sum_{l \in w} d_l g_l^p \neq 0$$

for some finite $w \subseteq [l^p, \omega)$ and $d_l \in \mathbb{Q}$. Hence

(2.4) $\quad \sum_{l \in w} d_l a_{jl} = 0$ for $j < k$ and $v = \sum_{l \in w} d_l b_l \neq 0$.

Multiplying this homogeneous system of equations and the inequality by a large enough natural number we may assume that

$$d_l \in n^p \mathbb{Z} \text{ for all } l \in w.$$

We now want to define $q > p$ with $q \in D_m^4$ and distinguish two cases. If $\sum_{l < l^p} b_l a_l^p \neq m$ then choose $n^q$ large enough such that $n^p | n^q$ and $\sum_{l < l^p} b_l a^p - m \not\equiv 0 \mod n^q$ and put $u^p = u^q, M^p = M^q, A^p = A^q$. Then $p < q$ and $\sum_{l < l^q} b_l a_l^q \not\equiv m \mod n^q$ hence $q \in D_m^4$ is a contradiction, see (2.2). If $\sum_{l < l^p} b_l a_l^p = m$, then choose $l^q > \sup(w \cup \{l^p\})$ and define $q$ such that

$$a_l^q(t) = \begin{cases} a_l^p(t) & \text{if } t \in [0, l^p) \\ d_l & \text{if } l \in w \\ 0 & \text{if } l \in [l^p, \omega) \setminus w. \end{cases}$$

Set $u^q = u^p \subseteq H_2$ and using (2.4) let $n^q$ be large enough such that $n^p | n^q$ but $v \not\equiv 0 \mod n^q$. It follows $p < q$ and

$$\sum_{l < l^q} b_l a_l^q = \sum_{l < l^p} b_l a_l^q + \sum_{l \in w} b_l d_l = m + v.$$

Hence $q \in D_m^4$ is another contradiction, see (2.2). The linear dependence (2.3) between the $h_l^p$'s and $g_l^p$'s is shown. Now we want to use (2.3) to derive a final contradiction for (2.2). For each finite $w \subseteq \omega$ we have a $\mathbb{Q}$-vector space $V_w = \langle h_l^p : l \in w \rangle$ of finite dimension $\leq k$. Hence there is an $r \in \omega$ and a finite $w^* \subseteq [l^p, \omega)$ such that $h_l^p$ ($l \in w^*$) is a maximal independent set - and $V_{w^*}$ has maximal dimension $|w^*| = r \leq k$. If $w^* \subseteq w \subseteq [l^p, \omega)$ for some finite $w$, then the sub-matrix $(\mathbf{H}_w) = (h_l^p, l \in w)$ of $(\mathbf{H})$ has finite column rank $r$, hence row rank $r$ as well and there is a subset $z \subset \{1, \ldots, k-1\}$ of size $r$ such that

$$\{\mathbf{a_j} \restriction w \; : \; j \in z\} \text{ is maximal independent.}$$

By (2.3) $\mathbf{b} \restriction w$ is a linear combination of the $\{\mathbf{a_j} \restriction w \; : \; j \in z\}$ and there are *unique* elements $c_l \in \mathbb{Q}$, $l \in z$ such that $\mathbf{b} \restriction w = \sum_{l \in z} c_l \mathbf{a_l} \restriction w$. If we increase $w$ we have the same coefficients by maximal independence. Hence

$$(2.5) \qquad \mathbf{b} \restriction [l^p, \omega) = \sum_{l \in z} c_l \mathbf{a_l} \restriction [l^p, \omega).$$

We can choose $m' \in \mathbb{N}$ large enough such that $m'c_l \in s'\mathbb{Z}$ for all $l \in z$. If $t \in W$ is large enough, then $m' | a_{lt}$ for all $l \in z$. Using (2.5) we get

$$b_t = \sum_{l \in z} c_l a_{lt} \in s'\mathbb{Z}$$

contradicting $W$. Hence $D_m^4$ is dense in $\mathfrak{F}$.

In order to show density of the last collection of subsets, suppose there are $\mathbf{d} \in H_1, t \in \mathbb{Z}$ and $n_0 \in \mathbb{N}$ such that

$$(2.6) \qquad D_{\mathbf{d}tn_0}^5 \text{ is not dense in } \mathfrak{F}.$$

Hence there is $p \in \mathfrak{F}$ such that

$$(2.7) \qquad \text{no } q \in \; D_{\mathbf{d}tn_0}^5 \text{ satisfies } p \leq q.$$

Let $u^p = \{\mathbf{c}^i = \sum_{j \in \omega} c_j^i \mathbf{e}_j : i < k\}$ and $l^p < l < \omega$. We want to consider extensions $p \leq q$ with $l^q = l$ and hence let

$$F_l = \{(y_{l^p}, \ldots, y_{l-1}) \in \mathbb{Z}^{l-l_p} : \sum_{j=l^p}^{l-1} c_j^i y_j = 0, i < k\}$$

which is a non-trivial subgroup of the free group $\mathbb{Z}^{l-l_p}$ for any large enough $l$. Also let

$$s(y_{l^p}, \ldots, y_{l-1}) = n_0 (\sum_{i < l^p} a_i^p \mathbf{b}^i + \sum_{i=l^p}^{l-1} y_i \mathbf{b}^i) - t(\sum_{i<l^p} a_i^p \mathbf{e}_i + \sum_{i=l^p}^{l-1} y_i \mathbf{e}_i) - \mathbf{d}.$$

We claim that

$$(2.8) \qquad (y_{l^p}, \ldots, y_{l-1}) \in F_l \Rightarrow s(y_{l^p}, \ldots, y_{l-1}) = 0 \text{ holds in } \mathbb{D}$$

If $s(y_{l^p}, \ldots, y_{l-1}) \neq 0$ for some $(y_{l^p}, \ldots, y_{l-1}) \in F_l$, then there is some $m \in \mathbb{N}$ such that

$$(2.9) \qquad s(y_{l^p}, \ldots, y_{l-1}) \not\equiv 0 \bmod m\mathbb{D}.$$

We now define some $q \in \mathfrak{F}$ taking

$$l^q = l, n^q = n^p \cdot m, u^q = u^p, M^q = \{m_{\mathbf{x}}^q = \sum_{i < l^q} x_i a_i^q : \mathbf{x} \in u^q\}$$

where

$$a_i^q = \begin{cases} a_i^p & \text{if} \quad i < l^p \\ y_i & \text{if} \quad l^p \leq i < l. \end{cases}$$

Clearly $q \in \mathfrak{F}$ and also $q \in D_{\mathbf{d}tn_0}$ from (2.9), hence $p \not\leq q$ from (2.7). On the other hand $\sum_{j=l^p}^{l-1} c_j^i a_j^q = 0$ from $F_l$ and definition of $a_i^q$ would imply $p \leq q$, a contradiction which proves the claim (2.8).

If we let
$$\mathbf{s}^i = \sum_{j \in \omega} s_j^i \mathbf{e}_j = n_0 \mathbf{b}^i - t\mathbf{e}_i \in \mathbb{D} \quad (l^p \leq i < \omega),$$
then the implication of (2.8) can be written as
$$\sum_{i=l^p}^{l-1} y_i \mathbf{s}^i = \mathbf{d} + t \sum_{i<l^p} a_i \mathbf{e}_i - n_0 \sum_{i<l^p} a_i^p \mathbf{b}^i.$$

From $(0, \ldots, 0) \in F_l$ follows

(2.10) $$n_0 \sum_{i<l^p} a_i^p \mathbf{b}^i = \mathbf{d} + t \sum_{i<l^p} a_i \mathbf{e}_i$$

and from

(2.11) $$(y_{l^p}, \ldots, y_{l-1}) \in F_l \text{ also follows } \sum_{i=l^p}^{l-1} y_i \mathbf{s}^i = 0.$$

If we view $\mathbf{s}^i = \sum_{j \in \omega} s_j^i \mathbf{e}_j$ as infinite row vector ($l^p \leq i < l$), then from the matrix
$$\begin{pmatrix} s_0^{l^p} & s_1^{l^p} & \cdots & s_k^{l^p} & \cdots \\ s_0^{l^p+1} & s_1^{l^p+1} & \cdots & s_k^{l^p+1} & \cdots \\ \vdots & & & & \\ s_0^{l-1} & s_1^{l-1} & \cdots & s_k^{l-1} & \cdots \end{pmatrix}$$
we have finite column vectors $\mathbf{s}_n = (s_n^i : l^p \leq i < l)$ for any $n \in \omega$. Let $\mathbf{c}^i \restriction [l^p, l)$ be the restriction of $\mathbf{c}^i$ viewed as an infinite column vector restricted to the coordinates $j$ such that $l^p \leq j < l$, then
$$\langle \mathbf{c}^i \restriction [l^p, l) : i < k \rangle$$
denotes the vector space over $\mathbb{Q}$ generated by these finite column vectors. We claim that
$$\mathbf{s}_n \in \langle \mathbf{c}^i \restriction [l^p, l) : i < k \rangle \quad \text{for all } n \in \omega.$$
Naturally $F_l \subseteq \mathbb{Z}^{l-l^p} \subseteq \mathbb{Q}^{l-l^p}$. If $\overline{F}_l = \langle F_l \rangle$ denotes the subspace of $\mathbb{Q}^{l-l^p}$ generated by $F_l$, then $\overline{F}_l = \langle \mathbf{c}^i \restriction [l^p, l) : i < k \rangle^\perp$ where orthogonality is defined naturally by
$$U^\perp = \{x \in \mathbb{Q}^{l-l^p} : x \cdot u = 0 \ \forall u \in U\}$$
for $U \subseteq \mathbb{Q}^{l-l^p}$ and the obvious scalar product $x \cdot u = \sum_{i \leq l-l^p} x_i u_i$. From (2.11) follows
$$\overline{F}_l = \langle \mathbf{s}_n : n \in \omega \rangle^\perp.$$
Using $\perp$ again, we have
$$\langle \mathbf{s}_n : n \in \omega \rangle^{\perp\perp} \subseteq \langle \mathbf{c}^i \restriction [l^p, l) : i < k \rangle^{\perp\perp}$$
which is
$$\langle \mathbf{s}_n : n \in \omega \rangle \subseteq \langle \mathbf{c}^i \restriction [l^p, l) : i < k \rangle$$
as $\dim \mathbb{Q}^{l-l^p}$ is finite. This shows the claim.

Now let $l$ be large enough such that $\langle \mathbf{c}^i \restriction [l^p, l) : i < k \rangle$ has maximal dimension $k' \leq k$ and let $\mathbf{c}^i \restriction [l^p, l)$ ($i < k'$) be a basis of this vector space. We now can write

$$\mathbf{s}_n = \sum_{i < k'} r_i^{nl} \mathbf{c}^i \restriction [l^p, l)$$

with *unique* coefficients $r_i^{nl} \in \mathbb{Q}$. By uniqueness these coefficients are independent of $l$ for any larger $l$, say that $r_i^{nl} = r_i^n$. In the system of equations

$$\mathbf{s}_n = \sum_{i < k'} r_i^n \mathbf{c}^i \restriction [l^p, l), \ (l^p \leq l < \omega, n \in \omega)$$

we can also eliminate $l$ and get

$$\mathbf{s}_n = \sum_{i < k'} r_i^n \mathbf{c}^i \restriction [l^p, \omega), \ n \in \omega.$$

From $\mathbf{s}^j$ and $\mathbf{b}^j = \sum_{n \in \omega} b_n^j \mathbf{e}_n$ we have that $s_n^j = n_0 b_n^j - t \delta_{jn} = \sum_{i < k'} r_i^n c_n^i$ for any $n \geq l^p$, hence $(n_0 \mathbf{b}^j - t \mathbf{e}_j) \restriction [l^p, \omega) \in \langle \mathbf{c}^i \restriction [l^p, \omega) : i < k' \rangle$ and

$$U = \langle n_0 \mathbf{b}^j - t \mathbf{e}_j : \ j \in \omega \rangle_* \subseteq \mathbb{D}$$

has finite rank. Hence $U$ is a free direct summand of $\mathbb{D}$, see Fuchs [8]. If $n_0$ does not divide $t$, then *modulo* $n_0 \mathbb{D}$ the image of $U$ is $\langle t \mathbf{e}_j + n_0 \mathbb{D} : \ j \in \omega \rangle_*$ and has infinite rank, which is impossible. Hence $n_0 | t$ and we rename $t n_0^{-1}$ by $t$. Using purity, we get that $U = \langle \mathbf{b}^j - t \mathbf{e}_j : \ j \in \omega \rangle_*$ is a free direct summand of $\mathbb{D}$ which contradicts our assumption that condition $(ii)(b)$ does not hold. Hence $D_{\mathbf{d} t n_0}$ is dense in $\mathfrak{F}$ indeed, see (2.6).

We are ready to apply Martin's axiom. There is a generic set $\mathbb{G} \subseteq \mathfrak{F}$ which meets the dense subsets of $\mathfrak{F}$ just constructed. We define $\mathbf{a} = \sum_{i \in \omega} a_i \mathbf{e}_i$ such that $a_i = a_i^p$ for any $p \in \mathbb{G}$ with $i < l^p$. Here we applied $D_{l_0}^3$ and note that $\mathbb{G}$ is directed, hence $\mathbf{a}$ is well-defined. Also $\mathbf{a} \in \mathbb{D}$ by $D_m^2$. Let $H_1' = \langle H_1, \mathbf{a} \rangle_* \subseteq \mathbb{D}$ be the *pure* subgroup of $\mathbb{D}$ generated by $H_1'' = H_1 + \mathbb{Z}\mathbf{a}$ and $\mathbb{H}' = (H_1', H_2)$. Then clearly $\mathbb{H} \subseteq \mathbb{H}'$ and we claim that $\mathbb{H}' \in \mathfrak{P}$. It is enough to show $(iii)$ for $\mathfrak{P}$. If $\mathbf{c} \in H_1''$ then $\mathbf{c} = k\mathbf{a} + \mathbf{e}$ for some $k \in \mathbb{N}, \mathbf{e} \in H_1$. If $\mathbf{y} \in H_2$, then consider $\Phi(\mathbf{c}, \mathbf{y}) = k\Phi(\mathbf{a}, \mathbf{y}) + \Phi(\mathbf{e}, \mathbf{y})$. From density of $D_\mathbf{y}^1$ and $p \in D_\mathbf{y}^1 \cap \mathbb{G}$ and the choice of $\mathbf{a}$ follows $\Phi(\mathbf{a}, \mathbf{y}) = m_\mathbf{y}^p \in \mathbb{Z}$ and therefore $\Phi(\mathbf{c}, \mathbf{y}) \in \mathbb{Z}$. The map $\Phi$ extends to $H_1'' \times H_2 \longrightarrow \mathbb{Z}$. If $\mathbf{x} \in H_1'$ then $t\mathbf{x} = \mathbf{h} \in H_1''$ for some $t \in \mathbb{N}$ and if $\mathbf{x} = \sum_{i \in \omega} x_i \mathbf{e}_i, \ \mathbf{h} = \sum_{i \in \omega} h_i \mathbf{e}_i$ then $t\mathbf{x} = \sum_{i \in \omega} t x_i \mathbf{e}_i = \sum_{i \in \omega} h_i \mathbf{e}_i$ and $h_i = t x_i$ for all $i \in \omega$. Hence

$$\Phi(\mathbf{h}, \mathbf{y}) = \Phi(t\mathbf{x}, \mathbf{y}) = \sum_{i \in \omega} t x_i y_i = t(\sum_{i \in \omega} x_i y_i) = t\Phi(\mathbf{x}, \mathbf{y}) \in t\widehat{\mathbb{Z}} \cap \mathbb{Z}$$

and by purity of $\mathbb{Z} \subseteq_* \widehat{\mathbb{Z}}$ also $t\Phi(\mathbf{x}, \mathbf{y}) \in t\mathbb{Z}$ and by torsion-freeness $\Phi(\mathbf{x}, \mathbf{y}) \in \mathbb{Z}$. We have seen that $\mathbb{H}' \in \mathfrak{P}$. Next we claim that

(2.12) by definition of $\mathbf{a}$ and $\mathbf{b}$ we have $z = \Phi(\mathbf{a}, \mathbf{b}) = \sum_{i \in \omega} b_i a_i \in \widehat{\mathbb{Z}} \setminus \mathbb{Z}$.

Note that $a_i \longrightarrow 0$ in the $\mathbb{Z}$-adic topology, hence $b_i a_i \longrightarrow 0$ and $z \in \widehat{\mathbb{Z}}$ is well-defined. If $z \in \mathbb{Z}$ and $n \in \mathbb{N}$ then $\sum_{i < k} b_i a_i \equiv z \mod n$ for any large enough $k$, which contradicts $D_{|z|}^4$.

Finally we show that $\sum_{i\in\omega} a_i\mathbf{b}^i \notin H_1'$. Otherwise there are $t, n \in \mathbb{N}$ and $\mathbf{d} \in H_1$ such that

(2.13) $$n\sum_{i\in\omega} a_i\mathbf{b}^i - t\mathbf{a} - \mathbf{d} = 0.$$

Let $p \in \mathbb{G} \cap D^5_{\mathbf{d}tn}$ from density of $D^5_{\mathbf{d}tn}$ and choose $m$ from the definition of $D^5_{\mathbf{d}tn}$. Hence

$$n\sum_{i<l^p} a_i^p\mathbf{b}^i - t\sum_{i<l^p} a_i^p\mathbf{e}_i - \mathbf{d} \in \mathbb{D} \setminus m\mathbb{D}.$$

On the other hand $a_i^p = a_i$ for all $i < l^p$ from $p \in \mathbb{G}$ and $m|n^p$ by $p \in D^5_{\mathbf{d}tn}$. The set $\mathbb{G}$ is directed, hence $m|a_i$ for all $i \geq l^p$. so $n\sum_{i\geq l^p} a_i\mathbf{b}^i \in m\mathbb{D}$ as well as $t\sum_{i\geq l^p} a^i\mathbf{e}_i \in m\mathbb{D}$. The last displayed expression becomes $n\sum_{i\in\omega} a_i\mathbf{b}^i - t\mathbf{a} - \mathbf{d} \in \mathbb{D} \setminus m\mathbb{D}$ which contradicts (2.13). The Main Lemma 2.2 is shown. $\square$

From the proof of the Main Lemma 2.2 we have an immediate

**Corollary 2.3.** *If $\mathbb{H} = (H_1, H_2) \in \mathfrak{P}$, $\mathbf{a} \in \mathbb{D}$ with $\Phi(\mathbf{a}, \mathbf{y}) \in \mathbb{Z}$ for all $\mathbf{y} \in H_2$ and $\mathbb{H}_1' = \langle H_1, \mathbf{a}\rangle_* \subseteq \mathbb{D}$ then $(H_1', H_2) \in \mathfrak{P}$, in particular $\Phi : H_1' \times H_2 \longrightarrow \mathbb{Z}$.*

In order to show Theorem 1.2 we want to use an *ad hoc* and preliminary definition. Here we also use that $\Phi$ is symmetric.

**Definition 2.4.** *A pair $\mathbb{H} = (H_1, H_2)$ of pure subgroups of $\mathbb{D}$ is a full pair if the following holds.*
  (i) *There is an increasing continuous chain $\mathbb{H}_\alpha = (H_{\alpha 1}, H_{\alpha 2}) \in \mathfrak{P}$ with $\alpha \in 2^{\aleph_0}$ and union $(H_1, H_2)$.*
  (ii) *If $\mathbf{b} \in P \setminus \mathbb{D}$ and $d \in \{1, 2\}$ there is $\mathbf{a} \in H_d$ such that $\Phi(\mathbf{a}, \mathbf{b}) \in \widehat{\mathbb{Z}} \setminus \mathbb{Z}$.*
  (iii) *If $\mathbf{b} \in \mathbb{D}$, then there is $d \in \{1, 2\}$ and $\mathbf{b} \in H_d$ or for some $\mathbf{a} \in H_{3-d}$ follows $\Phi(\mathbf{a}, \mathbf{b}) \in \widehat{\mathbb{Z}} \setminus \mathbb{Z}$.*
  (iv) *If $d \in \{1, 2\}$ and $\mathbf{b}^n \in H_d$, $(n \in \omega)$, there is $\mathbf{a} = \sum_{i\in\omega} a_i\mathbf{e}_i \in H_d$ such that*
      (a) *either $\sum a_i\mathbf{b}^i \in \mathbb{D} \setminus H_d$ or*
      (b) *or there is $t \in \mathbb{Z}$ such that $\langle \mathbf{b}^j - t\mathbf{e}_j : j \in \omega \rangle$ is a free direct summand of finite rank.*

.

**Lemma 2.5.** *(ZFC + MA) There is a full pair $\mathbb{H} = (H_1, H_2)$.*

**Proof.** Enumerate $P \setminus \mathbb{D} = \{\mathbf{b}_\alpha : \alpha \in 2^{\aleph_0}\}$ and $\mathbb{D}^\omega = \{(b_\alpha^n)_{n\in\omega} : \alpha \in 2^{\aleph_0}\}$ with $2^{\aleph_0}$ repetitions such that any element appears $2^{\aleph_0}$ times. We want to construct the $\mathfrak{P}$-chain inductively and let $(H_{01}, H_{02}) = (S, S)$. By continuity we only have to define $\mathbb{H}_{\alpha+1}$. Alternatively we switch between 1 and 2, say we are in case $H_{\alpha 1}$ and consider $\mathbf{b}_\alpha$ and $(\mathbf{b}_\alpha^n)_{n\in\omega}$. By the Main Lemma 2.2 there is $\mathbf{a}_\alpha \in \mathbb{D}$ such that $(H_{(\alpha+1)1}, H_{\alpha 2}) \in \mathfrak{P}$ where $H_{(\alpha+1)1} = \langle H_{\alpha 1}, \mathbf{a}_\alpha \rangle_* \subseteq \mathbb{D}$ and $\Phi(\mathbf{a}_\alpha, \mathbf{b}_\alpha) \in \widehat{\mathbb{Z}} \setminus \mathbb{Z}$. Moreover $(\mathbf{b}_\alpha^n)_{n\in\omega}$ satisfies condition $(ii)$ of the Main Lemma 2.2 for $\mathbf{b}^n = \mathbf{b}_\alpha^n$. A dual argument (case 2) provides $H_{(\alpha+1)2}$. Hence $\mathbb{H}_{\alpha+1} = (H_{(\alpha+1)1}, H_{(\alpha+1)2}) \in \mathfrak{P}$. This finishes the construction of $\mathbb{H}$ and Definition 2.4 is easily checked. $\square$

**Lemma 2.6.** *If $\varphi \in H_1^*$ for a full pair $\mathbb{H} = (H_1, H_2)$, then there is $\mathbf{b} \in H_2$ with $\varphi = \Phi(\ , \mathbf{b})$*

**Remark** A similar result holds for $\varphi \in H_2^*$.

**Proof.** Let $b_j = \mathbf{e}_j\varphi \in \mathbb{Z}$ for all $j \in \omega$, and set $\mathbf{b} = \sum_{j\in\omega} b_j\mathbf{e}_j \in P$. If $\mathbf{a} \in H_1 \subseteq \mathbb{D}$, then write $\mathbf{a} = \sum_{j\in\omega} a_j\mathbf{e}_j$ and by continuity $\mathbf{a}\varphi = (\sum_{j\in\omega} a_j\mathbf{e}_j)\varphi = \sum_{j\in\omega} a_j(\mathbf{e}_j\varphi) = \sum_{j\in\omega} a_j b_j = \Phi(\mathbf{a},\mathbf{b})$. Hence $\varphi = \Phi(\ ,\mathbf{b})$. If $\mathbf{b} \in P \setminus \mathbb{D}$, then by Definition 2.4 there is $\mathbf{x} \in H_1$ with $\mathbf{x}\varphi = \Phi(\mathbf{b},\mathbf{x}) \in \widehat{\mathbb{Z}} \setminus Z$ contradicting $\varphi \in H_1^*$, hence $\mathbf{b} \in \mathbb{D}$. Similarly by Definition 2.4 $(iii)$ we have $\mathbf{b} \in H_2$ and the lemma follows. □

The pair $\mathbb{H} = (H_1, H_2)$ in Lemma 2.6 satisfies conditions $(i)$ and $(iii)$ of the Theorem 1.2. Reflexivity follows easily as in [9] or [10] because the dual maps are induced by scalar multiplication. As a subgroup of $P$, each $H_i$ is $\aleph_1$-free, see Fuchs [8]. Slenderness can easily be checked and is left to the reader, hence $(ii)$ of Theorem 1.2 follows. Condition $(iv)$ was added for sake of completeness, it was the main goal in [9] or [10] and can be derived here using the arguments from there. The final condition $(v)$ will follow immediately from our next Lemma 2.7.

**Lemma 2.7.** *If $\mathbb{H} = (H_1, H_2)$ is a full pair and $\sigma \in \mathrm{End}\, H_1$ then there is $s \in \mathbb{Z}$ such that $\sigma - s1 \in \mathrm{Fin}\, H_1$ where $\mathrm{Fin}\, H_1$ is the ideal of $\mathrm{End}\, H_1$ of all endomorphisms of finite rank.*

**Proof.** If $\mathbf{e}_j\sigma = \mathbf{b}^j$, $j \in \omega$ then using that $\mathbb{H}$ is a full pair, we find $\mathbf{a} = \sum_{i\in\omega} a_i\mathbf{e}_i \in H_1$ such that Definition 2.4$(iv)$ holds. By by continuity also
$$\mathbf{b} = \mathbf{a}\sigma = (\sum_{i\in\omega} a_i\mathbf{e}_i)\sigma = \sum_{i\in\omega} a_n\mathbf{b}^i \in H_1$$
which shows that we are in case $(b)$ of Definition 2.4$(iv)$. The subgroup $U = \langle \mathbf{b}^j - t\mathbf{e}_j : j \in \omega \rangle$ is a free direct summand of finite rank of $\mathbb{D}$. However the image of $S = \bigoplus_{i\in\omega} \mathbf{e}_i\mathbb{Z}$ under $\sigma - t\,\mathrm{id}$ is in $U$, hence $S(\sigma - t\,\mathrm{id})$ has finite rank, and by continuity the same holds for $H_1(\sigma - t\,\mathrm{id})$, this is to say that $\sigma - t1 \in \mathrm{Fin}\, H_1$. □

## 3. Large reflexive groups

Let $\kappa$ be a fixed supercompact cardinal. Then there is a $\kappa$-complete (fine) ultra-filter $U$ over $\kappa$ such that the constant function
$$j : V \longrightarrow M = \mathrm{Ult}\,(V, U)\ (x \longrightarrow j(x))\ (j(x)_\alpha = x \text{ for all } \alpha \in \kappa)$$
is an elementary embedding of the universe $V$ into the ultra power $M$; for details see Kanamori [14, pp. 471, 298 - 306, 37 - 56]. If $\rho$ is a cardinal, then
$$\mathfrak{H}(\rho) = \{x \in V : |TC(x)| < \rho\}$$
is the set of all sets in $V$ hereditarily $< \rho$ where $TC(x)$ denotes the transitive closure of the set $x$.

**Theorem 3.1.** *If $\kappa$ is a supercompact cardinal and $H$ is a dual group of cardinality $\geq \kappa$, then for any $\chi < \kappa$ there is a direct summand $H'$ of $H$ with $\chi \leq |H'| < \kappa$.*

Then the following is immediate.

**Corollary 3.2.** *Every reflexive group of cardinality $\geq \kappa$, with $\kappa$ supercompact, has arbitrarily large summands $< \kappa$.*

**Proof of Theorem 3.1:** Let $H = G^* = \text{Hom}(G, \mathbb{Z})$ be as in the theorem. If $|G| = \lambda_1, |H| = \lambda_2$ then let $\lambda > 2^{\lambda_1 + \lambda_2}$ and assume $G = \lambda_1, H = \lambda_2$ as sets and $\chi < \kappa$. If $\mathfrak{P} = \mathfrak{P}_\kappa(\mathfrak{H}(\lambda))$ is the poset of all subsets of $\mathfrak{H}(\lambda)$ of cardinality $< \kappa$, then by the above there is a $\kappa$-complete (normal and fine) ultrafilter $D$ on $\mathfrak{P}$ with elementary embedding

$$(\mathfrak{H}(\lambda), \epsilon) \prec M := \text{Ult}(\mathfrak{P}, D).$$

From $H = G^*$ each $h \in H$ gives rise to a homomorphism

$$\Phi(h, \ ) : G \longrightarrow \mathbb{Z}$$

and $\Phi : H \oplus G \longrightarrow \mathbb{Z}$ is a bilinear form. Moreover

$$\Phi(h, \ ) = 0 \Rightarrow h = 0,$$

hence $\Phi$ is not degenerated. Let $\mathfrak{C}$ be the set of all $N \in \mathfrak{P}$ subject to the conditions

(i) $G, H, \Phi \in N$
(ii) $\chi + 1 \subseteq N$
(iii) $N$ is an elementary submodel of $(\mathfrak{H}(\lambda), \epsilon)$.
(iv) If $\tau = \text{otp}(N \cap \lambda)$ is the order type of $N \cap \lambda$, then $(N, \epsilon)$ is isomorphic to $(\mathfrak{H}(\tau), \epsilon))$, say by an isomorphism $j_N$.

By supercompactness $\mathfrak{C} \in D$ and by Łoś's theorem ([14, p. 47, Theorem 5.2] the desired properties of $\mathfrak{H}(\lambda)$ carry over to $N$. If $N \in \mathfrak{C}$ then define

$$H' = H \cap N \text{ and } G' = G \cap N.$$

From $\chi + 1 \subseteq N \in P$ and $\chi + 1 \subseteq \lambda_1 = G, \lambda + 1 \subseteq \lambda_2 = H$ follows $\chi + 1 \subseteq H'$ and $\chi + 1 \subseteq G'$, hence

$$\chi \leq |H'| < \kappa \text{ and } \chi \leq |G'| < \kappa$$

and by $(iii)$

(3.1) $$H' \subseteq H, \ G' \subseteq G \text{ are subgroups.}$$

Similarly, if $\Phi' = \Phi \restriction H' \oplus G'$, then

$$\Phi' : H' \oplus G' \longrightarrow \mathbb{Z}$$

and from $(iii)$ and $\Phi'$ we have

$$H' = G'^*.$$

We are ready to use an old trick from functional analysis to show that $H'$ is also a summand of $H$. Let

$$G'^\perp = \{h \in H : \Phi(h, G') = 0\} \text{ where } \Phi(h, G') = \{\Phi(h, g) : g \in G'\}.$$

Clearly $G'^\perp \subseteq H$, and consider any $h \in H' \cap G'^\perp$. We have $\Phi(h, G') = 0$ and from $h \in H'$ follows that in the submodel $N$ the following holds

$$N \models (\forall x \in G'^N \longrightarrow \Phi(h, x) = 0).$$

By $(iii)$ we also have

$$(\mathfrak{H}(\lambda), \epsilon) \models (\forall x \in G \longrightarrow \Phi(h, x) = 0),$$

hence $\Phi(h, \ ) = 0$ and $h = 0$ because $\Phi$ is not degenerated. We conclude

$$H' \cap G'^\perp = 0, G'^\perp \subseteq H.$$

In order to show

(3.2) $$H' + G'^\perp = H$$

we consider any $h \in H = G^*$ and let $\phi = \Phi(h,\ ) \restriction G'$ which belongs to $G'^*$. From (3.1) we find $h' \in H'$ such that $\Phi(h',\ ) = \phi$. If $g' \in G'$ we have

$$\Phi(h - h', g') = \Phi(h, g') - \Phi(h', g') = g'\phi - g'\phi = 0,$$

hence $h - h' \in G'^\perp$ and $h \in H' + G'^\perp$, and (3.2) follows. All together we see that $H'$ is a summand of $H$ of the right size. $\square$

Rüdiger Göbel
Fachbereich 6, Mathematik und Informatik
Universität Essen, 45117 Essen, Germany
e–mail: R.Goebel@Uni-Essen.De
and
Saharon Shelah
Department of Mathematics
Hebrew University, Jerusalem, Israel
and Rutgers University, Newbrunswick, NJ, U.S.A
e-mail: Shelah@math.huji.ae.il